\documentclass[12pt]{amsart}
\usepackage[english]{babel}
\usepackage[T1]{fontenc}

\usepackage{amsmath}
\usepackage{amsthm}
\usepackage{amssymb}
\usepackage[all]{xy}
\usepackage{mathrsfs}
\usepackage{enumerate}
\usepackage[notcite, final, notref]{showkeys}
\usepackage[dvips]{color}
\newtheorem{theorem}{Theorem}[section]

\newtheorem{lemma}[theorem]{Lemma}
\newtheorem{proposition}[theorem]{Proposition}
\newtheorem{corollary}[theorem]{Corollary}
\newtheorem{definition}[theorem]{Definition}

\theoremstyle{definition}
\newtheorem{remark}[theorem]{Remark}

\newcommand{\bn}{\bigskip\noindent}

\newcommand{\w}{\omega}

\usepackage{xcolor}

\newcommand{\cF}{\mathcal{F}}
\renewcommand{\bn}{\mathbb N}
\newcommand{\bc}{\mathbb C}
\newcommand{\ml}{ML_q(\mathbb{C})}
\newcommand{\br}{\mathbb R}

\newcommand{\llangle}{\langle \kern -0.2em \langle}	
\newcommand{\rrangle}{\rangle \kern -0.2em \rangle}
\newcommand{\inner}[1]{\left\langle  #1 \right\rangle }

\title[]{}

\author[N. Alpay]{Natanael Alpay}
\address{(NA) Schmid College of Science and Technology \\
Chapman University\\
One University Drive
Orange, California 92866\\
USA}
\email{nalpay@chapman.edu}

\author[K. Diki]{Kamal Diki}
\address{(KD) Schmid College of Science and Technology \\
Chapman University\\
One University Drive
Orange, California 92866\\
USA}
\email{diki@chapman.edu}

\title{On the Mittag Leffler Bargmann (MLB) transform}

\begin{document}
\maketitle

\begin{abstract}
We introduce the Segal-Bargmann transform associated to the Mittag Leffler Fock space and study how it will be connected to the Fourier transform.
We will discuss also the counterpart of the creation and annihilation operator in this setting using the Caputo and Liouville operators.
Finally, we give an extension of these results to the case of quaternions,
in particular in the slice hyperholomorphic setting.\\
\end{abstract}



\section{Introduction}

In 1961 Bargmann introduced in \cite{MR0157250} a Hilbert space of entire functions on which the creation and annihilation operators, namely $$\displaystyle M_zf(z):=zf(z)\quad \text{and} \quad Df(z):=\frac{d}{dz}f(z)$$ are closed, densely defined operators that are adjoints of each other and satisfy the classical commutation rule 
$$\left[D,M_z\right]=\mathcal{I}$$ 

where $\left[.,.\right]$ and $\mathcal{I}$ are respectively the commutator and the identity operator.
This space is called the Fock-Segal-Bargmann space. It consists of entire functions that are square integrable on the complex plane with respect to the normalized Gaussian measure.
It turns out that the standard Schr\"odinger Hilbert space on the real line is unitary equivalent to the Fock space via the so-called Segal-Bargmann transform, for more details we refer to the books \cite{qtm-book,  neretin, zhu}. 
A generalized  construction of the Fock space built upon a Mittag-Leffler reproducing kernel function is introduced in \cite{mlf}.
In this paper we will call this space a Mittag-Leffler-Fock (MLF) space.
The authors of \cite{mlf} obtained a new geometric description of this space which is based on a special weight extending the classical normalized Gaussian measure.
They also proved that classical complex monomials form an orthogonal basis and computed their norms with respect to the MLF space.
A general approach was investigated in \cite{ACK}.
In the present work, based on these recent ideas we are proposing a new Segal-Bargmann type transform associated to the MLF space.\\
j
To this end, our approach is built upon a new generating function involving normalized Hermite functions multiplied by a specific normalized basis of the Mittag-Leffler-Fock space.
In \cite{mlf} the counterpart of creation and annihilation operators were suggested to be related to the notion of Caputo derivatives. We study further the commutation rules of such operators in the present work. \\

On the other hand, an interesting theory of functions of a quaternionic variable was introduced in 2007 in \cite{struppa}.
These are the so called slice-hyperolomorphic functions, which can be considered as a generalization of  classical holomorphic functions to the quaternionic setting.
In recent years this new theory was extensively developed and found several applications in different directions including Schur analysis, quaternionic operator theory, and quaternionic quantum mechanics see \cite{ACS2016,MTS2017,ColomboSabadiniStruppa2011} .
The topic of Fock-Segal-Bargmann space and associated integral transforms in the slice hyperholomorphic case  were studied in  \cite{AlpayColomboSabadini2014,DG2017} (see also \cite{KA12,KA} for the slice polyanalytic case).
In this paper we extend these results based on a quaternionic type Mittag Leffler kernel.\\

The paper has the following structure: in the next section we briefly collect some basic facts and preliminary results on the Bargmann transform,  Mittag Leffler Fock (MLF) space, and the Caputo fractional derivative.
In section 3, we introduce ML-Bargmann kernel, and define the associated Mittag Leffler Bargmann (MLB) integral transform associated to the MLF. We compute the action of this transform on the Hermite functions, leading to an orthonormal basis of the MLF space. Moreover, we prove that the MLB transform is a unitary operator from the classical $L^2$ Hilbert space on the real line onto the MLF space, and compute its inverse.
Section 4 is devoted to study how the MLB transform is connected to the classical Fourier transform.
In section 5 we prove some commutation rules involving the creation and annihilation operators in this frame work, which are related to the Caputo fractional derivative.
Finally, in section 6, we give the counterpart of these results in the quaternionic slice hyper-holomorphic setting.

\section{Preliminaries}
In this section we revise the important mathematical tools that will be used in the sequel, in particular we revise the notion of Fock spaces and Segal Bargmann transforms
\cite{ alpaybook,MR0157250,zhu}
Then we recall the definition of the Mittag Leffler Fock space introduced in \cite{mlf}.
\begin{definition}[Fock space]\label{def:fock}
The Fock (or Bargmann-Fock-Segal) space $\mathcal F$ consists of entire functions $f$ that are square integrable with respect to the Gaussian measure, i.e. 
\begin{equation}
  \label{gauss123}
\mathcal{F} = \left\{ f\in \mathcal{H}(\mathbb{C}):\frac{1}{\pi}\iint_{\mathbb C}|f(z)|^2e^{-|z|^2}dA(z)<\infty\right\},
\end{equation}
where $dA(z)=dxdy$ is the Lebesgue measure with $z=x+iy$.
\end{definition}

The Fock space is a reproducing kernel Hilbert space with reproducing kernel
\begin{equation}
  \label{fock}
 K(z,w) =  e^{z\overline{\w}},\quad z,w\in \mathbb{C} .
\end{equation}
We remark that up to a positive multiplicative factor, $\mathcal F$ is the unique Hilbert space of entire functions
in which 
\begin{equation}
\label{fock-equation}
  \partial_z^*=M_z,
\end{equation}
where $\partial_z$ denote the derivative with respect to $z$.

\begin{definition}[Inner product on $L^2(\mathbb{R})$ and $\mathcal{F}$]
We consider the following Hilbert spaces
\begin{enumerate}[(i)]
\item The space $L^2(\mathbb{R})$ of square integrable functions with inner product
\begin{equation}\label{Eq:InnerProduct_l}
\inner{\phi, \psi}_{L^2(\mathbb{R})} = \int_{\br} \phi(t)\overline{\psi(t)}dt,\quad \psi,\phi\in L^2(\br).
 \end{equation}
\item The Fock space $\cF$ endowed with the inner product \begin{equation}\label{Eq:InnerProduct_Fock}
\inner{f, g}_{\mathcal{F}} = \frac{1}{\pi} \int_{\mathbb{C}} f(t)\overline{g(z)}  e^{-|z|^2} dxdy,\quad f,g\in\cF.
 \end{equation}
\end{enumerate}
\end{definition}


\begin{definition}[Bargmann transform]

The Bargmann transform ${\mathcal{B}}: L^2(\mathbb{R}) \to \mathcal{F}$ is defined by,
\begin{equation}
    \label{eq:b}
    {\mathcal{B}}[\phi](z)=\int_\mathbb{R}  \left( \sum_{n=0}^{\infty} \psi_n(x)\frac{z^n}{\sqrt{n!}} \right)  \phi(x)  dx,\quad \phi\in L^2(\br),
\end{equation}
where $\psi_n(x)$ are the normalized Hermite polynomials.
\end{definition}

\begin{definition}[MLF space]\label{def:mlf}
Let $q>0$, the Mittag Leffler Fock space (MLF space) is introduced in \cite{mlf} as follows
 $$\displaystyle ML_q(\mathbb{C})=\left\lbrace f\in\mathcal{H}(\mathbb{C}),\quad \frac{1}{q\pi}\int_{\mathbb{C}}|f(z)|^2|z|^{\frac{2}{q}-2}e^{-\frac{|z|^2}{q}}dA(z)< \infty \right\rbrace.$$
\end{definition}

It is importnat to note that the reproducing kernel function associated to $ML_q(\mathbb{C})$ can be expressed in terms of the classical Mittag Leffler function. In particular, this leads to the kernel 
\begin{equation}
K_q(z,w)=E_q(z\overline{w}):=\displaystyle \sum_{n=0}^\infty \frac{(z\overline{w})^n}{\Gamma(qn+1)},
\end{equation}
for every $z,w\in \mathbb{C}$.
We have also the following sequential characterization
$$\displaystyle ML_q(\mathbb{C})=\left\lbrace f(z)=\sum_{n=0}^{\infty}a_nz^n,\quad \sum_{n=0}^{\infty}\Gamma(qn+1) |a_n|^2< \infty \right\rbrace.$$
\begin{remark}
The particular case $q=1$ corresponds to the classical Fock space.
\end{remark}

Finally we revise from \cite{mlf,mlf-book} the concept of Caputo fractional derivatives denoted by $D^q_*$
for $q>0$.

\begin{definition}[Caputo fractional derivative]
For $n\in\bn$, $f\in\mathcal{C}^{(n)}(\br)$, the Caputo fractional derivative of order $q$, where $n-1<q\leq n$ is given by,
\[
D^q_* f(t) := 
\begin{cases}
\frac{1}{\Gamma(n-q)}\int_{0}^{t} \frac{f^{(n)}(x)}{(t-x)^{q+1-n}}dx, & n-1<q<n\\
f^{(n)}(t),& q=n\in\bn
\end{cases}
\]
\end{definition}

Using \cite[Theorem 5.1, Definition 8]{mlf} we obtain the equivalent definition for the Caputo fractional derivatives in the ML spaces that will be useful in the present work. 

\begin{definition}[Caputo fractional derivative in ML spaces]
Let $q>0$ and {$f\in ML^2_q(\bc \setminus \br)$ }
be given by the power series expansion $f(z) =\sum_{n=0}^{\infty}a_n z^{qn}$, then the Caputo fractional derivative of $f$ of order $q$ denoted by $D^{q}_{*}$ for $z\in \bc\setminus \br_{-}$ is given by
\[ 
D^{q}_{*} f(z) = \sum_{n=1}^{\infty} a_n \frac{\Gamma(qn+1)}{\Gamma(q(n-1)+1)}z^{q(n-1)}.
\]
\end{definition}

We note that the Mittag-Leffler function $E_q(z)$ satisfies the eigenvalue problem for the Caputo fractional derivative.



\section{The Mittag Leffler Bargmann transform}
In this section we introduce the Mittag Leffler Bargmann (MLB) integral transform associated to the MLF space $ML_q$ and study its main properties. Indeed, let us first consider the following generating function
\begin{definition}[ML-Bargmann kernel]
Let $q\in \mathbb R$, $q>0$, we introduce the ML-Bargmann kernel defined by 
\begin{equation}
	\label{eq:aq}
A_q(z,x)=A_q^z (x):=\displaystyle \sum_{n=0}^\infty \frac{\overline{z}^n}{\sqrt{\Gamma(qn+1)}}\psi_n(x),
\end{equation}
where $\{\psi_n\}_{n\geq 0}$ are the normalized Hermite functions. 
\end{definition}
\begin{definition}[MLB transform]
For any $\phi\in L^2(\mathbb{R})$ we define the MLB integral transform to be
\begin{equation}
\label{eq:bq}
B_q(\phi)(z):=\displaystyle \int_{\mathbb{R}}\overline{A_q(z,x)}\phi(x)dx.
\end{equation}
\end{definition}
\begin{remark}
The particular case $q=1$ in \eqref{eq:bq} leads to the classical Segal-Bargmann transform in formula \eqref{eq:b}.
\end{remark}
\begin{proposition}
\label{prop:ker}
For any $z,w\in\mathbb{C}$, it holds that 
\begin{equation}
\langle A_q^w, A_q^z \rangle_{L^2(\mathbb{R})}=K_q(w,z).
\end{equation}
In particular, we have also 
\begin{equation}
||A_q^z||^2_{L^2(\mathbb R)}=E_q\left(|z|^2\right).
\end{equation}
\end{proposition}
\begin{proof}
	Let $z,w\in\mathbb{C}$, we have $$\langle A_q^z, A_q^w \rangle_{L^2(\mathbb{R})}=\displaystyle \int_{\mathbb{R}}A_q^z(x)\overline{A^w_q(x)}dx.$$
	So we get, $$\langle A^z_q, A_q^w \rangle_{L^2(\mathbb{R})}=\int_{\mathbb{R}}\left(\overline{\sum_{k=0}^\infty \frac{z^k}{\sqrt{\Gamma(kq+1)}}\psi_k(x)}\right) \left({\sum_{j=0}^\infty \frac{w^j}{\sqrt{\Gamma(jq+1)}}\psi_j(x)}\right)dx.$$
	
	Hence by the convergence of the double series above and using the fact that Hermite functions are real valued, we obtain,
	$$\displaystyle \langle A^z_q, A_q^w \rangle_{L^2(\mathbb{R})}=\sum_{k,j=0}^{\infty}\frac{\bar z^k{w}^j}{\sqrt{\Gamma(kq+1) \Gamma(jq+1)}}\int_\mathbb{R}\psi_k(x)\psi_j(x)dx.$$
	
	However, we know that the normalized Hermite functions $(\psi_k)_{k\geq 0}$ form an orthonormal basis of the space $L^2(\mathbb{R})$ so that we have, $$\displaystyle \int_\mathbb{R}\psi_k(x)\psi_j(x)dx=\delta_{k,j},\quad k,j\in \mathbb N$$ where $\delta_{k,j}$ is the Kronecker symbol. 
	Thus, we obtain $$\langle A_q^z,A_q^w \rangle_{L^2(\mathbb{R})}=\displaystyle \sum_{k=0}^{\infty} \frac{\bar z^k {w}^k}{\Gamma(qk+1)}= K_q(w,z),\quad z,w \in\mathbb{C}.$$
	Finally, we note that
	$$||A^z_q||_{ML_q(\mathbb{C})}^2=\langle A_q^z, A_q^z \rangle_{ML_q(\mathbb{C})}=E_q(|z|^2).$$
\end{proof}

\begin{lemma}
For any $m\geq 0$, it holds that 
\begin{equation}
B_q(\psi_m)(z)=\displaystyle \frac{z^m}{\sqrt{\Gamma(q m+1)}},
\end{equation}
 for all $z\in\mathbb{C}$.
Moreover, we have \begin{equation}
||B_q(\psi_m)||_{ML_q(\mathbb{C})}=1=||\psi_m||_{L^2(\mathbb{R})}.
\end{equation}
\end{lemma}
\begin{proof}
	For a fixed $m\geq 0$, using the formulas \eqref{eq:aq} and \eqref{eq:bq}, we write:
	\begin{align*}
		B_q(\psi_m)(z) &= \int\limits_{\mathbb R} \overline{ A_q(z,x)} \psi_m(x)dx \\
		&=\int\limits_{\mathbb R}\left( \sum_{j=0}^{\infty} \frac{z^j}{\sqrt{\Gamma(qj+1)}} \psi_j(x)\right)\psi_m(x) dx,\\
		 &=\sum_{j=0}^{\infty}\frac{z^j}{\sqrt{\Gamma(qj+1)}}\left(\int\limits_{\mathbb R}\psi_j(x)\psi_m(x)dx\right).\\
	\end{align*}
Since $\psi_m$ are the Hermite functions, we get the following simplification
\begin{equation}
	\label{eq:bq_psi}
	B_q(\psi_m)(z)
	=\sum_{j=0}^{\infty}\frac{z^j}{\sqrt{\Gamma(qj+1)}}(\delta_{j,m})
	=\frac{z^m}{\sqrt{\Gamma(qj+1)}} := e_{m,q}(z). 
\end{equation}
Using \cite{mlf}, we know that 
\begin{equation}
	\label{eq:psi_mq}
	\left\{ e_{m,q} (z)= \frac{z^m}{\sqrt{\Gamma(mq+1)}}\; ;\; m\geq 0 \right\}
\end{equation}
form an orthonormal basis (shortened as ONB) of $ML_q(\mathbb C)$. Therefore
$$\displaystyle ||B_q(\psi_m)||_{ML_q(\mathbb C)} = 1 = ||\psi_m||_{L^2(\mathbb R)}$$.
\end{proof}

\begin{lemma}
For any $\psi\in L^2(\mathbb{R})$, 
it holds that,
\begin{equation}
||B_q(\psi)||_{ML_q(\mathbb{C})}=||\psi||_{L^2(\mathbb{R})}.
\end{equation}
\end{lemma}
\begin{proof}
Let $\psi\in L^2(\mathbb R)$, 
we can expand $\psi$ using the Hermite functions $\{\psi_m\}_{m\geq 0}$ which from an ONB of $L^2(\mathbb R)$. Thus, there exists unique coefficients $\alpha_m \in \mathbb{C},m\geq 0$, such that the following holds:
\begin{enumerate}[(i)]
	\item $\displaystyle \psi(x) = \sum_{m=0}^{\infty} \psi_m(x)\alpha_m$.
	\item $\displaystyle ||\psi||^2 = \sum_{m=0}^{\infty} |\alpha_m|^2$.
\end{enumerate}
Using the expansion above and formula given by \eqref{eq:bq}, we have
\begin{align*}
	B_q(\psi)(z)&=\int\limits_{\mathbb R} \overline{A_q(z,x)}\psi(x)dx\\
	&=\int\limits_{\mathbb R}\overline{A_q(z,x)}\left(\sum_{m=0}^{\infty} \psi_m(x)\alpha_m \right)dx\\
	&=\sum_{m=0}^{\infty}\left( \int\limits_{\mathbb R} \overline{A_q(z,x)}\psi_m(x)dx  \right) \alpha_m\\
	&=\sum_{m=0}^{\infty} \left( B_q (\psi_m)(z) \right) \alpha_m \\
	&=\sum_{m=0}^{\infty} e_{m,q}(z)\alpha_m.
\end{align*}
The last equality above is obtained using \eqref{eq:bq_psi}. 
However $\{ e_{m,q}\}_{m\geq 0}$ is ONB of $ML_q(\mathbb C)$, see
\cite{mlf}. Thus we will have 
$$ ||B_q(\psi)||_{ML_q(\mathbb C)}^2 = \sum_{m=0}^{\infty} |\alpha_m|^2 = ||\psi||^2_{L^2(\mathbb R)}.$$

Hence $\displaystyle B_q:L^2(\mathbb R) \to ML_q(\mathbb C)$ is an isometric operator.
\end{proof}

\begin{theorem}
The MLB transform $B_q$ is an isometric isomorphism mapping $L^2(\mathbb{R})$ onto the MLF space $ML_q(\mathbb{C})$.
\end{theorem}
\begin{proof}
In order to prove that $B_q$ is surjective, we apply the fact that for every $m\in \mathbb N$, we have $B_q(\psi_m)= e_{m,q}$. 
Indeed let $f\in ML_q(\mathbb C)$ such that 
$f(z)=\sum_{m=0}^{\infty}\alpha_m e_{m,q}(z) $
then we set
$ \psi = \sum_{m=0}^{\infty} \psi_m\alpha_m.$
It is clear that $\psi\in L^2(\mathbb R)$ and $f=B_q(\psi)$.
This shows that MLB transform 
$B_q:L^2(\mathbb R) \to ML_q(\mathbb C)$
is an isometric surjective (isomorphism) operator. 
\end{proof}

Let $q>0$, we will study the inverse and adjoint operators associated to the MLB transform $B_q$.  To this end, we introduce the transform $T_q:ML_q(\mathbb{C})\longrightarrow L^2$ defined by 
\begin{equation}
\label{eq:tf}
\displaystyle T(f)(x)=\int_{\mathbb{C}}{A_q(z,x)}f(z)|z|^{\frac{2}{q}-2}e^{-\frac{|z|^2}{q}}dA(z),
\end{equation}
for every $f\in ML_q(\mathbb{C})$ and $x\in\mathbb{R}$. In particular, we prove the following 
\begin{theorem}
For every $q\in\mathbb R$, $q> 0$, the MLB transform $B_q$ is an unitary operator whose inverse and adjoint operators are given by 
\begin{equation}
(B_q)^{-1}=(B_q)^{*}=T.
\end{equation}
\end{theorem}
\begin{proof}
Let $f\in ML_q(\mathbb C)$, since $B_q$ is an unitary operator it is enough to prove,
\begin{equation}
\label{eq:inv}
B_q\circ T(f) = f .
\end{equation}
In particular this would show that $B_q\circ T = I$.
We set $\phi=T(f)$, then using formulas \eqref{eq:bq} and \eqref{eq:tf} combined with Fubini's Theorem we have,
\begin{align*}
	B_q(\phi)(z) & = \int\limits_{\mathbb R} \overline{A_q(z,x)}\phi(x)dx\\
	&=\int\limits_{\mathbb R} \overline{A_q(z,x)}T(f)(x)dx\\
	&=\int\limits_{\mathbb R} \overline{A_q(z,x)} \left[  \int\limits_{\mathbb C} {A_q(w,x)} f(w)|w|^{\frac{2}{q}-2} e^{-\frac{|w|^2}{q}}dA(w)  \right]dx\\
	&=\int_{\mathbb{C}} f(w)|w|^{\frac{2}{q}-2} e^{-\frac{|w|^2}{q}}\left[ \int\limits_{\mathbb R}\overline{A_q(z,x)}{A_q(w,x)}  dx\right] dA(w).
\end{align*}

From the definition of the ML-Bargmann kernel and Proposition \ref{prop:ker} we know,
$$K_q(z,w) = \langle A_q^z,A_q^w\rangle_{L^2(\mathbb R)}=  \int\limits_{\mathbb R} A_q(z,x)\overline{A_q(w,x)}  dx,$$
and since $\overline{K_q(z,w)} = K_q(w,z)$. We have,
\begin{align*}
B_q(\phi)(z) 
&=\int_{\mathbb{C}} f(w)|w|^{\frac{2}{q}-2} e^{-\frac{|w|^2}{q}} \overline{K_q(w,z)}  dA(w)
\end{align*}	
Since 
$$\int_{\mathbb{C}} f(w)\overline{ K_{q,z}(w) } |w|^{\frac{2}{q}-2} e^{-\frac{|w|^2}{q}}dA(w) = \langle f, K_{q,z}\rangle_{ML_q(\mathbb{C})} = f(z),$$
we get 
$$B_q(\phi)(z) = f(z).$$
Therefore
$$(B_q\circ T)(f)(z) = f(z) \implies B_q\circ T = I.$$

Finally, using the fact that $B_q$ is an unitary operator, we have,
$$T=B_q^{-1}=B_q^{*}.$$

\end{proof}

\section{MLB and Fourier transform}

We start by recalling the Fourier transform of a signal $f:\mathbb R \to \mathbb C$ in $L^2(\mathbb R, dx)$ defined by
\[ \mathcal F (f)(\lambda) = \frac{1}{\sqrt{2\pi}} \int\limits_{\mathbb R} e^{-i\lambda x} f(x) dx. \]

\begin{lemma}\label{lemma:4.2}
	Let $q\in \mathbb \mathbb{R}$, $q > 0$. Then for every $z\in \mathbb C$ and $\lambda \in \mathbb R$, the action of the Fourier transform on the ML=Bargmann kernel is given by 
	$$ \mathcal{F}(A_q^{{z}})(\lambda) = A_q(i{z},\lambda). $$
	
\end{lemma}
\begin{proof}
	Applying the Fourier transform on \eqref{eq:aq} recalled by,
	\begin{equation*}
		A_q(z,x)=A_q^z (x):=\displaystyle \sum_{n=0}^\infty \frac{\bar z^n}{\sqrt{\Gamma(qn+1)}}\psi_n(x).
	\end{equation*}

We note that the function $A_q^{{z}}$ belongs to $L^2(\br)$ (see Proposition \ref{prop:ker}) so that we can apply the Fourier transform $\cF$. 
Moreover, using the well known fact $\cF(\psi_n)= (-i)^n \psi_n $,
we get
\begin{align*}
	\cF(A_q^{{z}})(\lambda) &=\sum_{n=0}^{\infty}\frac{\bar{z}^n}{\sqrt{\Gamma(qn+1)}}\cF(\psi_n)(\lambda)\\
	&=\sum_{n=0}^{\infty}\frac{\overline{z}^n}{\sqrt{\Gamma(qn+1)}}(-i)^n\psi_n(\lambda)\\
	&=\sum_{n=0}^{\infty}\frac{\overline{(iz)}^n}{\sqrt{\Gamma(qn+1)}}\\
	&={A_q(iz,\lambda)}.
\end{align*}

\end{proof}

\begin{theorem}\label{tm:4.3}
	For $q\in \mathbb R$, $q\geq 0$, the following diagram is commutative
$$S_q: \xymatrix{
	ML_q(\mathbb{C}) \ar[r] \ar[d]_{B_{q}^{-1}} & ML_q(\mathbb{C}) \\ L^2(\mathbb{R}) \ar[r]_{\mathcal{F}} &L^2(\mathbb{R}) \ar[u]_{B_q}
}$$ 
	so that the composition operator below is well defined
	\begin{equation}
		\label{eq:s}
		S_q = B_q \circ \mathcal F \circ B_q^{-1}.
	\end{equation}
	More precisely, for every $f\in ML_q(\mathbb C)$ it holds that
	$$ S_q(f)(z) = -f(-iz),\quad \forall z\in\mathbb C. $$ 
	
\end{theorem}
\begin{proof}
    Let $f\in \ml$.
	The proof have three steps. First we set 
	$$\phi(x)=B^{-1}_q(f)(x)=T(f)(x),$$
	then
	\begin{equation*}
		\phi(x)=\int\limits_{\mathbb C}{A_q(z,x)}f(z)|z|^{\frac{2}{q}-2}e^{-\frac{|z|^2}{q}}dA(z).
	\end{equation*}
	
	Secondly, we apply the Fourier transform on $\phi$ and get
	\begin{align*}
		\cF(\phi)(\lambda)&=\frac{1}{\sqrt{2\pi}}\int\limits_{\mathbb R}e^{-i\lambda x}\phi(x)dx\\
		&=\frac{1}{\sqrt{2\pi}}\int\limits_{\mathbb R}e^{-i\lambda x} \left( \int\limits_{\mathbb C} {A_q^z(x)}f(z)|z|^{\frac{2}{q}-2}e^{-\frac{|z|^2}{q}} dA(z) \right)dx\\
		&=\frac{1}{\sqrt{2\pi}}\int\limits_{\mathbb C}\left(\underbrace{\int\limits_{\mathbb R}e^{-i\lambda x}{A_q^z(x)}dx}_{:=I_q(z,\lambda)} \right)f(z)|z|^{\frac{2}{q}-2}e^{-\frac{|z|^2}{q}}dA(z).
	\end{align*}
	Using Lemma \ref{lemma:4.2} we have 
	\begin{align*}
	  \frac{1}{\sqrt{2\pi}}\; I_q (z,\lambda) &=\cF(A_q^{{z}})(\lambda)={A_q(iz,\lambda)}.
	\end{align*}
	So
	\begin{align*}
		\cF(\phi)(\lambda)&=\frac{1}{\sqrt{2\pi}}\int\limits_{\mathbb{C}}I_q(z,\lambda) f(z)|z|^{\frac{2}{q}-2}e^{-\frac{|z|^2}{q}}dA(z)\\
		&=\int\limits_{\mathbb{C}}{A_q(iz,\lambda)} f(z)|z|^{\frac{2}{q}-2}e^{-\frac{|z|^2}{q}}dA(z).
	\end{align*}
	Making the change of variable $z=-iw$ and $g(w)=-f(-iw)$ we get
	\begin{align*}
	\cF(\phi)(\lambda)&=\int\limits_{\mathbb{C}}{A_q(w,\lambda)} g(w)|-iw|^{\frac{2}{q}-2}e^{-\frac{|w|^2}{q}}dA(w)=B_q^{-1}(g)(\lambda).
	\end{align*}	
	Finally we apply $B_q$ and get
	$$ S_q(f)(z) = -f(-iz),\quad \forall z\in\mathbb C. $$ 

\end{proof}

\begin{corollary}
The Fourier transform is equivalent to the simple operator $S_q:ML_q(\mathbb C)\to ML_q(\mathbb C)$ under $MLB$ transform.
In particular we have
$$B_q^{-1} \circ S_q \circ B_q = \mathcal F.$$
\end{corollary}
\begin{proof}
	We note that by Theorem \ref{tm:4.3} we have 
	$$S_q = B_q \circ \mathcal F \circ B_q^{-1}.$$
	Then, if we multiply by the $B_q^{-1}$ transform from the left and by $B_q$ transform from the right we obtain 
	$$\mathcal{F}=B_q^{-1} \circ S_q \circ B_q$$ 
\end{proof}

\section{Creation and annihilation operators}

In this section we study the notion of creation and annihilation operators on the MLF space based on the MLB transform.
We start by introducing the notion of Caputo derivatives,
\begin{definition}
For $f(z)=\sum_{n=0}^{\infty}z^{qn}a_n$ we have as Caputo Derivative
\begin{equation}
    D^q_*(f)(z)=\sum_{n=1}^{\infty}a_n\frac{\Gamma(qn+1)}{\Gamma(q(n-1)+1)} z^{q(n-1)}.
\end{equation}
\end{definition}

and the multiplication operator $M_{z^q}:=z^q$.
We note that Caputo derivative $D_*^q$, and $M_{z^q}$ could be used to define an analog for the creation and annihilation operator. Indeed we will compute the following commutator operator,
\begin{equation}
\label{eq:com}
\begin{aligned}
[D_*^q,M_{z^q}]&=D^q_*M_{z^q}(f)(z)-M_{z^q}D_*^q(f)(z)\\
&=D^q_*(z^qf)(z)-z^qD_*^q(f)(z).
\end{aligned}
\end{equation}

\begin{proposition}
For every $f(z)=\sum_{n=0}^{\infty}z^{qn}a_n$, it holds that
\begin{align*}
[D_*^q,M_{z^q}](f)(z)
=a_0\Gamma(q+1)+\sum_{n=1}^{\infty}a_n\beta_{m,q}\;z^{qn},
\end{align*}
with the coefficients
$$\beta_{m,q}=
\left[ {\frac{\Gamma(q(n+1)+1)}{\Gamma(qn+1)}-\frac{\Gamma(qn+1)}{\Gamma(q(n-1)+1)}}
\right].
$$
\end{proposition}
\begin{proof}
For $f(z)=\sum_{n=0}^{\infty}a_nz^{qn}$, set $g(z)=z^qf(z)$, then
\[
g(z)=\sum_{m=0}^{\infty}a_mz^{q(m+1)}=\sum_{n=1}^{\infty}a_{n-1}z^{qn}=\sum_{n=0}^{\infty}b_nz^{qn},\quad b_n=
\begin{cases}
a_{n-1}&, n\geq 1,\\
0&, \text{otherwise}.
\end{cases}
\]

So the Caputo derivative applied to $g(z)$ would be,
\begin{align*}
D_*^q(g(z)) &= D_*^q(z^qf)(z)\\
&=\sum_{n=1}^{\infty} b_n\frac{\Gamma(qn+1)}{\Gamma(q(n-1)+1}z^{q(n-1)},\\
&=\sum_{n=1}^{\infty}a_{m-1}\frac{\Gamma(qn+1)}{\Gamma(q(n-1)+1)}z^{q(n-1)}.
\end{align*}

We also have,
\begin{align*}
z^q D_*^q(f)(z) =\sum_{n=1}^{\infty} a_n\frac{\Gamma(qn+1)}{\Gamma(q(n-1)+1}z^{qm}.
\end{align*}
Therefore the commutator \eqref{eq:com} would be
\begin{align*}
[D_*^q,M_{z^q}](f)
&=D^q_*(z^qf)(z)-z^qD_*^q(f)(z)\\
&=a_0\Gamma(q+1)+\sum_{n=1}^{\infty}z^{qn}a_n\left[ \underbrace{\frac{\Gamma(q(n+1)+1)}{\Gamma(qn+1)}-\frac{\Gamma(qn+1)}{\Gamma(q(n-1)+1)}}_{:=\beta_{m,q}}
\right]\\
&=a_0\Gamma(q+1)+\sum_{n=1}^{\infty}a_n\beta_{m,q}\;z^{qn}.
\end{align*}
\end{proof}

\begin{remark}
The commutation rule proved in the previous result can be seen as an extension of the classical commutation rule on the Fock space.\\
\end{remark}


\begin{proposition}
For different values of $q$ we have the following expressions:
\begin{enumerate}
    \item 
    When $q=1$, then $\beta_{m,1}=1$ which gives us back the identity and the classical case.\\
    \item
    When $q=2$, then $\beta_{m,2}=2(4n+1)$ which results in the commutator
    \[ [D_*^2,M_{z^2}](f)(z)  = 2\left(f(z^2)+2z\frac{d}{dz}f(z^2)\right). \]
    \item 
    When $q=3$, then $\beta_{m,3}=27(3n+1)+6$ which results in the commutator
    \[ [D_*^3,M_{z^3}](f)(z)  = 3\left(2f(z^3)+6z\frac{d}{dz}f(z^3)+3z^2\frac{d^2}{dz^2}f(z^3)\right). \] 
    \item When $q=4$, we have,
 \[
 [D^{*}_{q},M_{z^4}](f)(z)=24f(z^4)+120z\frac{d}{dz}(f(z^4))+56z^2\frac{d^2}{dz^2}(f(z^4))+16 z^3\frac{d^3}{dz^3}(f(z^4)).\\
 \]
 \end{enumerate}
\end{proposition}

\begin{proof}
 The calculation for $q=1,2$ is easy to verify, we provide the calculation for $q=3$ and $q=4$ follows a similar pattern.
 
 We first remark the following
 \[ z^2\frac{d}{dz^2}(z^{3m})=z^2\frac{d}{dz}(3mz^{3m-1})
 = (3m)^2z^{3m}-3mz^{3m}.\]
 So 
\begin{equation}
\label{eq:q3}
(3m)^2z^{3m} = z^2\frac{d}{dz^2}(z^{3m}) 3mz^{3m}. 
\end{equation}
 
 For $q=3$, we have the following where the last line is obtained by \eqref{eq:q3}
 \begin{align*}
    [D^{*}_{q},M_{z^3}](f)(z)&= a_0q!+\sum_{m=1}^{\infty}z^{qm} a_m\left[ \frac{(q(m+1))!}{(qm)!}-\frac{(qm)!}{(q(m-1))!} \right]\\
    &= a_0 6 + \sum_{m=1}^{\infty} z^{3m}a_m \left[ \frac{(3(m+1))!}{(3m)!}-\frac{(3m)!}{(3(m-1))!} \right]\\
    &= a_0 6 + \sum_{m=1}^{\infty} z^{3m}a_m (27m(3m+1)+6)\\
    &= 6\left(a_0 + \sum_{m=1}^{\infty} z^3a_m \right) + 27 \sum_{m=0}^{\infty} a_m m(3m+1) z^{3m}\\
    &= 6\left(\sum_{m=0}^{\infty} z^3a_m \right) + 9 \sum_{m=0}^{\infty} a_m (3m)  z^{3m} + 9\sum_{m=0}^{\infty} a_m (9m^2)  z^{3m}\\
    &= 6\left(\sum_{m=0}^{\infty} z^3a_m \right) + 9 \sum_{m=0}^{\infty} a_m (3m)  z^{3m} + 9\sum_{m=0}^{\infty} a_m \left((3m)z^{3m} + z^2\frac{d^2}{dz^2}(z^3n)\right)\\
    &=6f(z^3)+18z\frac{d}{dz}(f(z^3))+9z^2\frac{d^2}{dz^2}(f(z^3)),
\end{align*}
by \eqref{eq:q3}.

\end{proof}


    \begin{remark}
    When $q\in\mathbb N,q>0$, then
    $$
    [D_*^q,M_{z^2}] (f)(z) = \sum_{k=0}^{q-1} C_k z^{k}\frac{d^{k}}{dz^k}f(z^q).
    $$

We plan to prove this conjecture in a forthcoming work.

    \end{remark}

\section{Quaternionic MLB transform}

In this section we adapt the notion of MLF and MLB transform to the case of quaternions, in particular in the slice hyperholomorphic setting. 

The noncommutative field $\mathbb{H}$ of quaternions consists of elements of the form
$p=x_0  + x_{1} i +x_{2} j + x_{3} k$, $x_i\in \mathbb{R}$, $i=0,1, 2,3$,
where the imaginary units $i, j, k$ satisfy the relations
$$i^2=j^2=k^2=-1, \ ij=-ji=k, \ jk=-kj=i, \ ki=-ik=j.$$ The real number $x_0$ is called real part of $p$
while $ x_{1} i +x_{2} j + x_{3} k$ is called imaginary part of $p$.
We define the norm of a quaternion $p$ as
$|p|=\sqrt{x_0^2+x_1^2 +x_2^2+x_3^2}$.
Note that if $I\in \mathbb{S}$, then $I^{2}=-1$.  For any fixed $I\in\mathbb{S}$ we define $\mathbb{C}_I:=\{x+Iy; \ |\ x,y\in\mathbb{R}\}$, which can be identified with a complex plane.

Any non real quaternion $p$ is uniquely associated to the element $I_p\in\mathbb{S}$
defined by $I_p:=( i x_{1} + j x_{2} + k x_{3})/|  i x_{1} + j x_{2} + k x_{3}|$ and $p$ belongs to the complex plane $\mathbb{C}_{I_p}$.

Obviously, the real axis belongs to $\mathbb{C}_I$ for every $I\in\mathbb{S}$ and thus a real quaternion can be associated to any imaginary unit $I$.
Moreover, we have $\mathbb{H}=\cup_{I\in\mathbb{S}} \mathbb{C}_I$.

In this part of the paper we are interested in the so called slice regular functions of a quaternion variable which are defined below. 
For more information on these functions theory and for their applications, we refer to \cite{ColomboSabadiniStruppa2011}, \cite{gss}.

\begin{definition}\label{slice regular}
	Let $U$ be an open set of
	$\mathbb{H}$. A real differentiable
	function $f:U \to \mathbb{H}$ is said to be (left) slice regular if,
	for every $I \in \mathbb{S}$, its restriction $f_I$ to $U_I=U\cap\mathbb C_I$ satisfies
	$$
	\overline{\partial}_If(x+Iy)=\frac{1}{2}\Big(\frac{\partial}{\partial x}
	+I\frac{\partial}{\partial y}\Big)f_I(x+Iy)=0.
	$$
	The set of slice regular functions on $U$ will be denoted by $\mathcal{SR}(U)$. If $U=\mathbb{H}$ then $\mathcal{SR}(\mathbb{H})$ will be called as the space of entire slice regular functions.
\end{definition}

In \cite{AlpayColomboSabadini2014,DG2017} the notions of Fock spaces and Segal Bargmann transform were studied in the setting of slice hyperholomorphic functions of quaternionic variable.

\begin{definition}[Quaternionic Fock space and Segal Bargmann transform]
Let $I$ an imaginary unit in $\mathbb{S}$, the Quaternionic Fock space is defined by
$$\displaystyle \mathcal{F}_{slice}(\mathbb{H})=\left\lbrace f\in\mathcal{SR}(\mathbb{H}),\quad \frac{1}{\pi}\int_{\mathbb{C}_I}|f(p)|^2| e^{-|p|^2}dA_{I}(p)< \infty \right\rbrace,$$
where $dA_I(p)=dxdy $ when $p=x+Iy$.\\

The Quaternionic Bargmann transform ${\mathcal{B}_{slice}}: L^2(\mathbb{R},\mathbb{H}) \to \mathcal{F}_{slice}(\mathbb{H})$ is defined by,
\begin{equation}
    \label{eq:b}
    {\mathcal{B}_{slice}}[\phi](p)=\int_\mathbb{R}  \left( \sum_{n=0}^{\infty} \psi_n(x)\frac{p^n}{\sqrt{n!}} \right)  \phi(x)  dx,\quad \phi\in L^2(\br,\mathbb{H}),
\end{equation}
where $\psi_n(x)$ are the normalized Hermite polynomials.

\end{definition}

\begin{remark}
It is important to note that the space $\mathcal{F}_{slice}(\mathbb{H})$ does not depend on the choice of the imaginary unit $I$ (see \cite{AlpayColomboSabadini2014}).
Moreover, in \cite{DG2017}, it was proved that the Quaternionic Bargmann transform $\mathcal{B}_{slice}$ is a surjective isometry from $L^2(\br,\mathbb{H})$ onto $\mathcal{F}_{slice}(\mathbb{H})$.  
\end{remark}



In this section we will extend the MLF spaces and MLB transform studied in previous sections to the case of quaternions.

\begin{definition}[QMLF space]\label{def:mlf}
Let $q\in \mathbb R$, $q>0$, and $I$ an imaginary unit in $\mathbb{S}$. Then, the Quaternionic slice hyperholomorphic Mittag Leffler Fock space (QMLF space) is defined by
$$\displaystyle ML_q(\mathbb{H})=\left\lbrace f\in\mathcal{SR}(\mathbb{H}),\quad \frac{1}{q\pi}\int_{\mathbb{C}_I}|f(p)|^2|p|^{\frac{2}{q}-2}e^{-\frac{|p|^2}{q}}dA(p)< \infty \right\rbrace.$$
\end{definition}

The reproducing kernel function associated to $ML_q(\mathbb{H})$ can be expressed in terms of the non-commutative Mittag Leffler function. In particular, this leads to the kernel 
{
\begin{equation}
K_q(p,s)=E_q^*(p\overline{s}):=\displaystyle \sum_{n=0}^\infty \frac{p^n\overline{s}^n}{\Gamma(qn+1)},
\end{equation}}
for every $p,s\in \mathbb{H}$.
\begin{remark}
Using standard arguments we can justify that the QMLF space $ML_q(\mathbb{H})$ does not depend on the choice of the slice $\mathbb{C}_I$.
\end{remark}

\begin{definition}[QML-Bargmann kernel]
Let $q\in \mathbb R$, $q>0$, we introduce the QML-Bargmann kernel defined by 
\begin{equation}
	\label{eq:q_aq}
\mathcal{A}_q(p,x)=\mathcal{A}_q^p (x):=\displaystyle \sum_{n=0}^\infty \frac{\overline{p}^n}{\sqrt{\Gamma(qn+1)}}\psi_n(x), \quad p\in \mathbb{H}, x\in \mathbb{R}
\end{equation}
where $\{\psi_n\}_{n\geq 0}$ are the normalized Hermite functions. 
\end{definition}

We consider the space $L^2(\mathbb{R},\mathbb{H})$ consisting of functions $\phi:\mathbb {R}\to \mathbb{H}$ such that
\[ \int_{\mathbb{R}} |\phi(x)|^2dx <\infty. \]

On this space the inner is given by
\[
\inner{\phi,\psi}_{L^2(\mathbb{R})}= \int_{\mathbb{R}} \overline{\psi(x)}\phi(x)dx.
\]
\begin{definition}[QMLB transform]
For any $\phi\in L^2(\mathbb{R},\mathbb{H})$ we define the QMLB integral transform to be
\begin{equation}
\label{eq:q_bq}
\mathcal{B}_q(\phi)(p):=\displaystyle \int_{\mathbb{R}}A_q(p,x)\phi(x)dx.
\end{equation}
\end{definition}

\begin{proposition}
\label{prop:q_ker}
For any $p,s\in\mathbb{H}$, it holds that 
{\begin{equation}
\langle \mathcal{A}_q^s, \mathcal{A}_q^p \rangle_{L^2(\mathbb{R,\mathbb{H}})}=K_q(p,s).
\end{equation}}
In particular, we have also 
\begin{equation}
||\mathcal{A}_q^p||^2_{L^2(\mathbb R,\mathbb{H})}=E_q\left(|p|^2\right).
\end{equation}
\end{proposition}

\begin{lemma}
For any $m\geq 0$, it holds that 
\begin{equation}
\mathcal{B}_q(\psi_m)(p)=\displaystyle \frac{p^m}{\sqrt{\Gamma(q m+1)}},
\end{equation}
 for all $p\in\mathbb{H}$.
Moreover, we have \begin{equation}
||\mathcal{B}_q(\psi_m)||_{QML_q(\mathbb{H})}=1=||\psi_m||_{L^2(\mathbb{R,,\mathbb{H}})}.
\end{equation}
\end{lemma}

\begin{theorem}
The QMLB transform $\mathcal{B}_q$ is an isometric isomorphism mapping $L^2(\mathbb{R},\mathbb{H})$ onto the QMLF space $QML_q(\mathbb{H})$.
\end{theorem}

\begin{remark}
It would be interesting to investigate the counterpart of the creation and annihilation operators for quaternionic MLF spaces.
To this end it would be important to develop a notion of the quaternioinc Caputo derivatives.
\end{remark}



\begin{thebibliography}{}

\bibitem{alpaybook}
Alpay, Daniel. \textit{An Advanced Complex Analysis Problem Book: Topological Vector Spaces, Functional Analysis, and Hilbert Spaces of Analytic Functions.} Birkhäuser, 2015.


\bibitem{ACS2016} Alpay D., Colombo F., Sabadini I.,
       {\it Slice Hyperholomorphic Schur Analysis}.
     Operator Theory. Advances and Applications, Vol. 256. Birkhäuser, Basel Mathematics (2016)



\bibitem{AlpayColomboSabadini2014} Alpay D., Colombo F., Sabadini I., Salomon G., {\it The Fock space in the slice hyperholomorphic Setting}. In Hypercomplex Analysis: New perspectives and applications. Trends Math. pp. 43-59. (2014)

\bibitem{ACK}
Alpay Natanael, Paula Cerejeiras, and Uwe Kähler. \textit{"Generalized Fock space and fractional derivatives with Applications to Uniqueness of Sampling and Interpolation Sets."} arXiv preprint arXiv:2112.07883 (2021).

\bibitem{MR0157250}
V.~Bargmann.
\newblock \textit{On a {H}ilbert space of analytic functions and an associated integral
  transform.}
\newblock {\em Comm. Pure Appl. Math.}, 14:187--214, 1961.

\bibitem{qtm-book}
Hall, Brian C. \textit{Quantum theory for mathematicians.} Vol. 267. New York: Springer, 2013.


\bibitem{MTS2017} Muraleetharan B., Thirulogasanthar K., Sabadini I., {\it A representation of Weyl–Heisenberg Lie algebra in the quaternionic setting}.
 Ann. Physics. 385, pp. 180-213 (2017)


\bibitem{ColomboSabadiniStruppa2011} Colombo F., Sabadini I., Struppa D.C., {\it Noncommutative Functional Calculus, Theory and Applications of Slice Hyperholomorphic Functions}. Progress in Mathematics. Vol. 289. Birkh\"auser, Basel. (2011)


\bibitem{CD2017} Cnudde L., De Bie B., {\it Slice Segal-Bargmann transform}. J. Phys. A. Vol. 50, No. 25 (2017)


\bibitem{KA12}
De Martino, Antonino, and Kamal Diki. \textit{"On the polyanalytic short-time Fourier transform in the quaternionic setting."} arXiv preprint arXiv:2110.04520 (2021).

\bibitem{KA}
De Martino, Antonino, and Kamal Diki. \textit{"On the Quaternionic Short-Time Fourier and Segal–Bargmann Transforms."} Mediterranean Journal of Mathematics 18.3 (2021): 1-22.


\bibitem{DG2017} Diki, K., Ghanmi, A., {\it A Quaternionic Analogue of the Segal-Bargmann Transform}.  Complex Anal. Oper. Theory. Vol. 11, pp. 457-473. (2017)


\bibitem{struppa}
Daniele C. Struppa and Gentili, Graziano . \textit{"A new theory of regular functions of a quaternionic variable."} Advances in Mathematics 216.1 (2007): 279-301.


\bibitem{mlf}
Rosenfeld, Joel A., Benjamin Russo, karl Warren E. Dixon. \textit{‘The Mittag Leffler reproducing kernel Hilbert spaces of entire and analytic functions’}. Journal of Mathematical Analysis and Applications 463.2 (2018): 576-592. Web.

\bibitem{gss} G. Gentili, C. Stoppato, D. C. Struppa, {\em Regular Functions of a Quaternionic Variable}, Springer
Monographs in Mathematics, Springer, Berlin-Heidelberg, 2013.

\bibitem{mlf-book}
Gorenflo, R. k.a\`. \textit{Mittag-Leffler Functions, Related Topics and Applications.} Springer Berlin Heidelberg, 2014. Web. Springer Monographs in Mathematics.


\bibitem{neretin}
Neretin, Yu A. \textit{Lectures on Gaussian integral operators and classical groups.} Vol. 13. European Mathematical Society, 2011.


\bibitem{zhu}
Zhu, K. (2012). Fock Spaces. In: Analysis on Fock Spaces. Graduate Texts in Mathematics, vol 263. Springer, Boston, MA. 


\end{thebibliography}
\end{document}